\documentclass[12pt,a4paper]{article}
\usepackage{amsmath}
\usepackage{amssymb}
\usepackage{graphicx}
\usepackage{color}
\usepackage{bbold}
\usepackage{dsfont}

\usepackage[draft=false]{hyperref}

\newcommand{\levy}{L\'{e}vy }
\renewcommand{\a}{\alpha}
\newcommand{\R}{{\mathbb R}}

\newcommand{\E}{{\rm E}}
\newcommand{\e}{{\rm E}}
\renewcommand{\P}{{\rm P}}
\newcommand{\p}{{\rm P}}
\newcommand{\D}{{\mathrm d}}
\newcommand{\RR}{\mathbb{R}}
\newcommand{\CC}{\mathbb{C}}
\newcommand{\DD}{\mathbb{D}}
\newcommand{\sDD}{\mathbb{S}}
\newcommand{\ZZ}{\mathbb{Z}}
\newcommand{\matI}{I}
\newcommand{\matO}{0}
\renewcommand{\sp}{\mathrm{sp}}

\newcommand{\indic}{\mathbb 1}
\newcommand{\1}[1]{\indic_{\{#1\}}}

\renewcommand{\u}{{\mbox{\tiny $+$}}}

\newcommand{\bs}{\boldsymbol}
\newcommand{\vect}[1]{\boldsymbol #1}
\newcommand{\Dpi}{\Delta_{\bf\pi}}

\newcommand{\vone}{\vect 1}
\newcommand{\vzero}{\vect 0}

\newcommand{\hR}{R}
\newcommand{\wF}{\widetilde F}

\newcommand{\wL}{\widetilde L}

\newcommand{\wtau}{\widetilde \tau}
\newcommand{\wA}{\widetilde A}

\newcommand{\wR}{\widetilde R}
\newcommand{\wx}{\widetilde X}

\newcommand{\whG}{\widehat G}
\newcommand{\whH}{\widehat{H}}

\newcommand{\whJ}{\widehat{J}}

\newcommand{\whR}{\widehat{R}}
\newcommand{\whX}{\widehat{X}}

\newcommand{\vligne}[1]{\begin{bmatrix} #1 \end{bmatrix}}

\newtheorem{defn}{Definition}[section]
\newtheorem{lem}[defn]{Lemma}
\newtheorem{prop}[defn]{Proposition}

\newtheorem{thm}[defn]{Theorem}
\newtheorem{cor}[defn]{Corollary}
\newtheorem{rem}[defn]{Remark}
\newtheorem{remark}[defn]{Remark}
\newtheorem{ass}[defn]{Assumption}

\newcommand{\qed}{\hfill $\square$}
\newenvironment{proof}{
      \noindent {\bf Proof }}{\qed
      \vspace{0.25\baselineskip}
}
\newcommand{\debproof}{\begin{proof}}
\newcommand{\finproof}{\end{proof}}

\newcommand{\ct}{{\bf The non-lattice case}}
\newcommand{\cta}{{\bf The non-lattice Case A}}
\newcommand{\ctb}{{\bf The non-lattice Case B}}
\newcommand{\dt}{{\bf The lattice case}}

\definecolor{darkmagenta}{rgb}{0.5,0,0.5}
\definecolor{darkgreen}{rgb}{0,0.6,0}
\definecolor{darkblue}{rgb}{0,0,0.6}
\definecolor{darkred}{rgb}{0.8,0,0}
\definecolor{mellow}{rgb}{.847, 0.72, 0.525}

\usepackage{pifont}

\usepackage[normalem]{ulem}

\begin{document}

\title{One-sided Markov additive processes with lattice and non-lattice increments}

\author{
Jevgenijs Ivanovs\thanks{Aarhus University, Department of Mathematics, \texttt{jevgenijs.ivanovs@math.au.dk}}
\and
Guy Latouche\thanks{Universit\'e libre de Bruxelles,
D\'epartement d'informatique, CP 212, Boulevard du Triomphe, 1050
Bruxelles, Belgium, \texttt{latouche@ulb.ac.be}.}
\and
Peter Taylor\thanks{University of Melbourne, School of Mathematics and Statistics, Vic 3010, \texttt{taylorpg@unimelb.edu.au}}
}
\date{\today}

\maketitle

\begin{abstract}
Dating from the work of Neuts in the 1980s, the field of
matrix-analytic methods has been developed to analyse discrete or
continuous-time Markov chains with a two-dimensional state space in
which the increment of a {\it level} variable is governed by an
auxiliary {\it phase} variable.  
More recently, matrix-analytic techniques have been applied to general
Markov additive models with a finite phase space.
The basic assumption underlying these developments is that the process
is skip-free (in the case of QBDs or fluid queues)
or that it is  {\it one-sided}, that is it is jump-free in one direction.

From the Markov additive perspective, traditional matrix-analytic
models can be viewed as special cases:  for M/G/1 and GI/M/1-type
Markov chains, increments in the level are constrained to be {\it
  lattice} random variables and for fluid queues, they have to be
piecewise linear.

In this paper we discuss one-sided lattice and non-lattice Markov
additive processes in parallel. Results that are standard in one
tradition are interpreted in the other, and new perspectives
emerge. In particular,  using three fundamental
matrices, we address hitting, two-sided exit, and
creeping probabilities.
%
  \end{abstract}


\section{Introduction}
\label{sec:intro}

A Markov additive process (MAP) is a bivariate Markov process
$(X_t,J_t), t\geq 0$ living on $\mathbb R\times E$ and adapted to some
right-continuous, complete filtration $\mathcal F_t$. The components
$X$ and $J$ are referred to as the \emph{level} and \emph{phase}
components respectively, and it is assumed throughout this work that
the phase space $E=\{1,\ldots,N\}$ is finite. The defining property of
a MAP states that, for any time $T$ and any phase~$i$, conditional on
$\{J_T=i\}$, the process $(X_{T+t}-X_T,J_{T+t}),t\geq 0$ is
independent of $\mathcal F_T$ and has the law of $(X_t-X_0,J_t),t\geq
0$ given \mbox{$\{J_0=i\}$}.  
Thus, the increments of the level process are governed by the phase
process $J$, with the latter evolving as a continuous-time Markov
chain with some transition rate matrix~$Q$; without real loss of
generality we assume that $J$ is irreducible.  

In general, a MAP can be thought of as a \emph{Markov modulated \levy
  process} in the sense that $X$ evolves as a \levy process $X^{(i)}$
while $J=i$, with additional jumps distributed as $U_{ij}$ when $J$
switches from $i$ to~$j \not= i$, all components assumed to be
independent.

If $X$ lives on a lattice, then $(X,J)$ is a bivariate Markov chain
and its  transition rate matrix has a repetitive block structure.  We
usually assume that $X$ takes integer values; we refer to these
processes as \emph{lattice} as opposed to \emph{non-lattice} MAPs.  
Avram and Vidmar \cite{avvi19} used analogous terminology when they derived scale functions for (scalar) \levy processes with integer-valued increments.

We
focus here on {\em skip-free downward} processes.  In the lattice case, this means that the process may jump from
level $n$ to levels $n-1$, $n$, $n+1$, $n+2$, \ldots but not to any
level lower than $n-1$.  Early references are
Neuts~\cite{neuts_book,neut89}, where skip-free upward Markov chains
are also extensively studied. If $X$ is allowed to take any real value, 
we assume that the process has no negative jumps.  In both the lattice
and non-lattice cases, we shall call such MAPs {\em one-sided}.

Our objective is to bring together results that have been obtained
from different angles in the literature, thereby emphasising the unity
resulting from the Markov additive point of view.  At the same time,
differences in the model definitions lead to differences in the
structure of the results, and these will be clarified as
well. 

In the next section, we shall introduce our main notation and discuss the
assumptions that will be made throughout the paper.  In Section
\ref{s:ladder} we consider the ladder-height process of first passages
to lower levels. This is the opportunity to introduce $G$, one of
the fundamental matrices used in the analysis of MAPs.  Section
\ref{s:occupation} is about occupation times, and that is where we
introduce another important matrix $H$. Then, in Section \ref{s:R},
we discuss returning to zero and the matrix $R$.  In Section
\ref{s:two-sided}, we analyse exit probabilities from a bounded
interval and discuss scale matrices.  We look at applications in
Section \ref{s:applications}: to creeping probabilities, the distribution of
extrema when the phase process terminates in finite time, and 
additional results for occupation times.  Section \ref{s:conclusion} is
devoted to some concluding observations.

\section{Preliminaries}
\label{s:prelim}

\subsection{Notation and assumptions.} \label{sec:defn}

Putting the Markov additive property to work requires repeated conditioning on the phase, which leads us to use a matrix formalism.

We write $\e[\cdot ; J_t]$ to denote the $N\times N$ matrix with $ij$-th element
\[
\e_{i}[\cdot;J_t=j]=\e[\,\cdot\,  \1{J_t=j}|J_0=i],
\]
similarly, $\p[\cdot ; J_t]$ stands for a matrix with corresponding probabilities.  Unless otherwise specified, we assume that $X$ starts at~0.  The identity and the zero matrices are denoted by $\matI$ and $\matO$ respectively, the (column) vectors of ones and zeros by $\vone$ and $\vzero$ and we write $\Delta_{\bs v}$ for a diagonal matrix with the vector $\bs v$ on the diagonal.

Because we wish to allow for {\it defective} processes, we allow $Q$
to be non-conservative, that is, its row sums may be strictly less
than zero, and the non-negative vector $\bs q=-Q\vone$ provides {\it
  killing rates} at which the process transitions into a {\it cemetery
  state $\dagger$} where it remains for all subsequent time. We usually avoid including $\dagger$ as a possible value of $J_t$, considering it implicit in the statement $\{J_t=i\}$ that the process has survived up to time~$t$. 

\subsection{The non-lattice case.} \label{sec:nonl}

In the non-lattice case we assume that $X_t \in \RR$ and that jumps
are positive. It is well known (see, for instance, Asmussen and Kella~\cite{ak00}) that such a MAP is characterized by a matrix-valued function 
\begin{align}\label{eq:F}
 F(\alpha)=Q\circ (\e [ e^{\alpha
  U_{ij}}])+\Delta_{\vligne{\psi_1(\alpha) & \ldots & \psi_N(\alpha)}},
\end{align}
where $\circ$ stands for entry-wise (Hadamard) matrix
multiplication,  $\psi_i(\alpha)\equiv\log \e [e^{\alpha X^{(i)}_1}]$ is
the Laplace exponent of $X^{(i)}$ and  $U_{ii}=0$. The function
$F(\alpha)$ is defined at least for $\Re(\a)\leq 0$, but may have a
larger domain. Specifically it is defined over the set ${\cal D}_F$
for which the jumps, both those of the L\'evy processes $X^{(i)}$ and
the $U_{ij}$s, have a finite moment generating function.
%

The L\'evy-Khintchine formula provides a general expression for
$\psi_i(\alpha)$:
\begin{equation}
   \label{eq:psi}
\psi_i(\alpha)=\frac{1}{2}\sigma_i^2\alpha^2+a_i\alpha+\int_0^\infty(e^{\alpha
    x}-1-\alpha x\1{x<1})\nu_i(\D x)\qquad \alpha\in {\cal D}_F,
\end{equation} 
where $(a_i,\sigma_i,\nu_i(\D x))$ is the \levy triplet corresponding to $X^{(i)}$, see \cite[Theorem~1.6]{kyprianou}.

With $F$ defined as in (\ref{eq:F}), we have
\begin{align}\label{eq:EXJ}
\e[e^{\alpha X_t};J_t]=e^{tF(\alpha)}
\end{align}
for all $\alpha\in {\cal D}_F$. 

\begin{ass}
   \label{a:nosub}
   Significant simplifications occur if we assume that none of the
   $X^{(i)}$s is a subordinator, that is, none of the $X^{(i)}$s have
   paths that are almost surely non-decreasing. This leads us to
   assume that for each $i$, either $\sigma^2_i > 0$, or $\sigma^2_i = 0$ and
   $a_i < 0$. 
\end{ass}

Finally, note that $F(0)=Q$, and that we can introduce additional killing with rates $\bs q'$ by substituting $F(\alpha)-\Delta_{\bs q'}$ for $F(\alpha)$, because $U_{ii}=0$ by convention.

\subsection{The lattice case.} \label{sec:latt}

For lattice MAPs, we assume that $X_t \in\mathbb Z$, and so $(X,J)$ is
a continuous-time Markov chain on the state space
$\mathbb Z \times \{1,\ldots,N\}$. The MAP property implies that state
transitions are translation-invariant in the level component, and we
let $A_m$ be the block of the transition rates from level $k$ to level
$k+m$.  As usual, the diagonal entries of $A_0$ are strictly negative
and are such that the row sums of the transition matrix are negative
or zero. The analogue of the non-lattice assumption that there are no
negative jumps is the assumption that the level component is
\emph{skip-free downwards}, with $A_m=\matO$ for $m=-2,-3,\ldots$.

From the viewpoint of a Markov modulated \levy process, our lattice case corresponds to the situation where all the \levy processes $X^{(i)}$ are compound Poisson processes with jumps $B_i$, and where all the $B_i$s and $U_{ij}$s take values in $\{-1,0,1,\ldots\}$. 

Even though the Laplace transform of the increments can still be written in the form~\eqref{eq:F} for the lattice case, it is more convenient to use a discrete generating function, defined on a set ${\cal D}_F$ that contains $\{ z:|z|\leq 1,z\neq 0\}$ by
\begin{align}\label{eq:F_lattice}
F(z)=z\sum_{m=-1}^{\infty} z^mA_m,\end{align}
the factor $z$ outside the sum ensuring that $F$ does not have a
singularity at~0. Then  
\begin{align}\label{eq:EXJ_lattice}
\e[z^{X_t};J_t]=e^{tF(z)/z}
\end{align}
for any $z \in {\cal D}_F$. 
The transition rate matrix of $J$ is $Q=F(1)=\sum_{m=-1}^{\infty}A_{m}$. As with the non-lattice case, one can introduce additional killing by writing $F(z)-z\Delta_{\bs q'}$ instead of $F(z)$.



In the lattice case, our assumption that $Q$ is irreducible does not rule out pathological examples such as discussed in Asmussen~\cite[p.\ 314]{APQ} where the MAP is constrained to move within a finite number of levels, or Latouche and Taylor~\cite{LT97}, where $(X,J)$ is made up of several non-communicating processes.  
\begin{ass}
   \label{cond:1} 
To avoid such processes,
we assume that  the Markov chain $(X,J)$ is irreducible.
\end{ass}


\subsection{Asymptotic drift}\label{sec:drift}

Consider, for the moment, that $\bs q = \vzero$, so that the process
is non-defective.  It is well-known (see, for example, \cite[Cor.\
XI.2.8]{APQ}) that there is a constant $\mu\in [-\infty,\infty]$, such
that $\lim_{t\rightarrow\infty}X_t/t = \mu$ a.s.\ for every initial
state; in particular, $\lim_{t\rightarrow\infty}X_t$ is $\infty$ or
$-\infty$ as $\mu>0$ or $\mu<0$.  If $\mu=0$, then
$\limsup_{t\rightarrow\infty} X_t=\infty$ and
$\liminf_{t\rightarrow\infty} X_t=-\infty$ a.s.\  under our
assumptions, see \cite[Prop.\ XI.2.10]{APQ}; we refer to this as the
{\it zero-drift} case.
In the sequel, we will call a MAP {\it transient} if $\bs q \not= \vzero$ or $\bs q = \vzero$ and $\mu \not = 0$, and {\it null recurrent} if $\bs q = \vzero$ and $\mu=0$.  Null recurrent
processes often require special attention.

We also mention that
$\mu=\e_{\bs\pi} [X_1]=\sum_i\pi_i\e[X_1|J_0=i]$, where $\bs\pi$ is the
stationary distribution of~$J$.  In the {non-lattice} case this
can be written in the form $\mu=\bs \pi F'(0)\bs 1$, and in the {lattice case}
$\mu=\bs \pi (\sum_{k=-1}^\infty k A_k)\bs 1$.

\section{The first passage Markov chain and the matrix $G$}
\label{s:ladder}

Among the objects of interest in the analysis of MAPs, we find
the first passage time into a level
\begin{equation} \label{eq:taux}
 \tau_x =\inf\{t\geq 0:X_t=x\}, \qquad x\in\mathbb R.
\end{equation}
The Markov additive property, together with the assumption that $X$ is one-sided, implies that the phase observed at the first passage times $\tau_{-x}$ into negative levels forms a Markov chain with time index $x$, which we shall call the \emph{first passage Markov chain}. In the lattice case, this is a discrete-time Markov chain characterised by  
\begin{align}
\label{eq:Gdef}
 \p[J_{\tau_{-k}}]=G^k, \qquad k=0,1,\ldots
\end{align}
for some \emph{transition probability} matrix $G$, while in the {non-lattice case} it is a continuous time Markov chain with transition probabilities of the form 
\begin{align}
\label{eq:Gdef_cont}
 \p[J_{\tau_{-x}}]=e^{Gx}, \qquad x\geq 0
\end{align}
for some \emph{transition rate} matrix~$G$. 
We will distinguish between the two regimes: 
\begin{equation}
   \label{e:regimes}
\begin{split}
{\rm Regime}\ C1: {} & \mbox{the process $(X_t,J_t)$ is non-defective and $\mu\leq 0$, and}
\\
{\rm Regime}\ C2: {} & \mbox{either the process $(X_t,J_t)$ is defective or $\mu>0$.}
\end{split}
\end{equation}
In Regime C1, $\tau_{-x}$, defined in (\ref{eq:taux}) is almost surely finite for any $x>0$, and the first passage Markov chain is non-defective, whereas in Regime C2 there is a non-zero probability that $\tau_{-x}=\infty$, for any positive $x$ and the first passage Markov chain is defective.


\subsection{The lattice case} \label{sec:latt1}

In addition to the continuous-time Markov chain $(X_t,J_t)$, it is useful here to consider the discrete-time Markov chain $(\widetilde X_n,\widetilde J_n)$ for $n=0,1,\ldots$ embedded at jump
  epochs,  discussed for example in Neuts~\cite[Section 1.7]{neuts_book} and in Ramaswami and Taylor ~\cite[Section 4]{rata96}.  Let $\Delta_A$ be the diagonal matrix with the same diagonal as $-A_0$. The transition probabilities $\widetilde A_m$ corresponding to an increase of $m$ levels for $\wx$ are given by
\begin{equation}
\label{eq:jumpchain}
\widetilde A_m=\Delta_A^{-1}A_m, \quad m \not=0, \qquad \qquad \wA_0=
\Delta_A^{-1}A_0+I. 
\end{equation}
By conditioning on $\wx_1$, one observes that 
\[   
 \widetilde A_{-1}+\widetilde A_0G+ \widetilde A_1G^2+\ldots=G.
\]   
Multiplying this on the left by $\Delta_A$ and using~\eqref{eq:jumpchain}, we see that
\begin{equation}\label{eq:G}
 A_{-1}+A_0G+A_1G^2+\ldots=\matO.
\end{equation}
This is the well-known equation for $G$ established by
Neuts~\cite{neut89}; In other words, $G$ is a \emph{right} solution to
the equation $F(\cdot)=\matO$.  In fact, $G$ is the
minimal non-negative solution to this equation, it is the
unique stochastic solution in Regime C1, and the unique substochastic,
but not stochastic, solution in Regime C2 \cite[Theorem 4.3]{blm05}.

Except in some very special cases, equation \eqref{eq:G} cannot be solved
analytically.   This  has paved the way to the development of very
efficient numerical algorithms. 
%
Software routines for carrying out the computations are downloadable
from Van Houdt's webpage~\cite{vanh14}.

\subsection{The non-lattice case} \label{sec:nonlatt1}
The matrix $G$ is a conservative generator in Regime C1 and has at least
one row sum that is strictly negative in Regime C2. 
Moreover, under our assumptions that $Q$ is irreducible and none of
the $X^{(i)}$s is a subordinator, $G$ is {irreducible} in both
cases. 

Let $(a_i,\sigma_i,\nu_i(\D x))$ be a \levy triplet corresponding to
$X^{(i)}$ chosen as in~\eqref{eq:psi}, and let $U(\D x)$ be the matrix
of distributions of~$U_{ij}$.  Here again, $G$ is a right
solution of $F(\cdot)=\matO$ in the sense that
\begin{equation}
   \label{eq:FG}
   \frac{1}{2}\Delta^2_{\bs\sigma}    G^2+\Delta_{\bs a}G+\int_0^\infty\Delta_{\bs \nu(\D
     x)}\left(e^{Gx}-\matI-Gx\1{x<1}\right)+\int_0^\infty (Q\circ U(\D
   x)) e^{G x}=\matO,
\end{equation}
see, for example, \cite[Eq.\ (6)]{breuer_lambda} or~\cite[Thm. 2]{lambda}.

As a solution of~\eqref{eq:FG}, the matrix $G$  is identified in the
following way.  Define $\mathcal M_1$ to be the set of $N\times N$
matrices with a simple eigenvalue at 0 and all others having negative
real parts, and $\mathcal M_2$ to be the set of matrices of which all
eigenvalues have negative real parts. Via a Jordan decomposition, we
observe that the integrals in~\eqref{eq:FG} converge absolutely for
any matrix in $\mathcal M_1$ or $\mathcal M_2$. Then by~\cite[Thm. 2]{lambda} in Regime C1, $G$ is the unique matrix in
$\mathcal M_1$ solving~\eqref{eq:FG}, and in Regime C2, it is the
unique matrix in $\mathcal M_2$ solving~\eqref{eq:FG}.  Finally, the
matrix $G$ can be obtained numerically using iterative methods based
on equation~\eqref{eq:FG}, see~\cite{asmu95, breuer_lambda, matth17, blm21}, or via
the spectral method developed in~\cite{lambda}.  The latter, however,
may not be practical in the case when the number $N$ of phases is large.

\begin{rem} \em
Interestingly, equation (\ref{eq:FG}) is much older than is widely
believed. It was stated in general form by Ezhov
and Skorokhod as far back as 1969, see~\cite[Eq.\ (30)]{ezhov_skorokhodII}. 
Their equation, though, has a typographical error and the authors did not elaborate on its derivation.
Among the possible methods for deriving it, we mention (a) the
infinitesimal generator approach~\cite{prabhu_storage,rogers,
  breuer_lambda}, (b) factorization of L\'evy processes and
embedding~\cite{asmu95, dieker_mandjes}, and (c) the spectral method
based on martingale calculations~\cite{lambda}.
\end{rem}



\section{The expected occupation time and the matrix $H$}
\label{s:occupation}

\subsection{Occupation time at a level}\label{sec:L}
The occupation time at a level is an important concept. It is easy to understand in the lattice case.  In the non-lattice case it was defined in~\cite{ivanovs_palmowski} based on the corresponding theory for \levy processes.  Moreover, this concept laid the basis for the results in~\cite{ivanovs_potential}.

In the non-lattice case, we define the occupation time at level $x$ on
the interval $[0,t)$ to be the vector whose $j$th component is 
\begin{equation}
\label{eq:ldefunbounded}
 L_{j}(x,t)=\lim_{\epsilon\downarrow 0}\frac{1}{2\epsilon}\int_0^t\1{|X_s-x|<\epsilon,J_s=j}\D s, 
\end{equation}
where the limit is taken uniformly on compacts in $L^2$, implying that
$L_{j}(x,t)$ is continuous in~$t$.  If $X^{(j)}$ is a process with
bounded variation on compacts, then $L_j(x,t)$ may be written in terms
of the number of crossings of level $x$ in phase $j$:
\begin{equation}\label{eq:Ldefbounded}
L_{j}(x,t)=\frac{1}{|a_j|}\#\{s\in[0,t):X_s=x,J_s=j\}.
\end{equation} 
Under Assumption \ref{a:nosub}, the linear drift coefficient $a_j$ of
such a phase must be negative, since $\sigma_j^2=0$.

In the lattice case, we use
\begin{equation}
   \label{e:time}
L_{j}(x,t) =\int_0^t\1{X_s=x,J_s=j}\,\D s.
\end{equation}
which is the time spent by $(X,J)$ in state $(x,j)$ up to time $t$.

In both the lattice and non-lattice cases, $L_{j}(x,t),t\geq 0$ is an $\mathcal F_t$-adapted process, that increases only when $(X,J)=(x,j)$, and which has an important additive property.
Consider a stopping time $\tau$ such that $X_\tau=y$ and $J_\tau=k$
for some~$y$ and $k$.  Conditional on the event $J_\tau=k$, the
shifted process $L_{j}(x,\tau+s)-L_{j}(x,\tau),s\geq 0$ is independent
of $\mathcal F_\tau$ and has the law of $L_{j}(x-y,s),s\geq 0$ given
$J_0=k$. 

Furthermore, for any measurable non-negative function $f$
which is such that $f(\cdot,\dagger)=0$, the identity 
\begin{equation}\label{eq:occ_density}
 \int_0^t f(X_s,J_s)\D s=\sum_{j\in E} \int_{\mathbb R}f(x,j)L_{j}(x,t)\D x 
\end{equation}
holds almost surely, with the integral on the right-hand side replaced by $\sum_{x\in\mathbb Z}$ in the lattice case.

\subsection{Definition and properties}
The expected occupation time at a level $x\in\R$ on the infinite time horizon is defined by a matrix $H(x)$ whose $(i,j)$th component is
\[H(x)_{ij}=\e_i [L_j(x,\infty)].\]
We shall write $H(x) = \e[L(x,\infty)]$ and use the simplified notation $H=H(0)$. 

Our first result is that $H(x)$ has finite entries for all $x$ unless $(X,J)$ is null-recurrent, in which case $H(x)$ has infinite entries for all $x$. 
\begin{lem}
   \label{t:Hfinite}
   {If the process is transient, either because $\mu\not = 0$ or
     because the generator $Q$ is defective, then} the matrix $H(x)$
     has finite entries for all $x \in \mathbb R$.  {If the
       process is null-recurrent, that is, if the matrix $Q$ is
       non-defective and $\mu = 0$, then $H(x)$ is infinite for all
       $x$.}
\end{lem}
\begin{proof}
The defective case is trivial, and so we work with the non-defective case. 
Assume that $J_0=i$ for some fixed phase $i$, 
and define
\begin{align*}
t'_1  & = \inf_{t\geq 1}\{t: (X_t,J_t)=(0,i)\}   \\
t'_{n+1} &  = \inf_{t\geq t'_n + 1}\{t: (X_t,J_t)=(0,i)\} \qquad
            \mbox{for $n\geq 1$,}
\end{align*}
with $t'_{n}=\infty$ if these events do not occur.  Also, define
$N_i = \sup\{n:t'_n < \infty\}$; this is the number of times that
$(X_t,J_t)$ returns to $(0,i)$ after at least a
unit time interval.

If $\mu \not=0$,  it is known that $|X_t|\rightarrow\infty$ a.s. as $t\rightarrow
\infty$, and that $X$ almost surely leaves the neighbourhood
of~$0$ in finite time, never to return.  Therefore, 
$N_i$ is finite a.s., the
conditional independent increment property of the process ensures that
$N_i$ is geometrically-distributed and, by the additivity of $L$, we write
\[
\e_i[L_{i}(0,\infty)]=\e_i[L_{i}(0,1)] \e_i[N_i] < \infty.
\]
This further implies that all entries of $H$ are finite.  We also
immediately have 
\begin{equation}\label{eq:HxH}
H(x)=\p[J_{\tau_x}]H, \qquad x\in\mathbb R,
\end{equation}
which connects the matrix-valued function $H(x)$ with the basic first
passage time into a level and proves that $H(x)$ is finite for all $x \in \mathbb R$.

If the process is null-recurrent, then $N_i=\infty$ a.s.,  $H$ is
infinite, and so is $H(x)$ for all $x$.
\end{proof}

For $x\leq 0$, note that in the non-lattice case 
\begin{equation}
   \label{e:HandG}
H(x)=e^{-Gx}H
\end{equation}
and so it is continuous in this domain. Moreover, $H(x)$ is continuous for $x>0$, because a.s. none of
the underlying \levy processes $X^{(i)}$ can jump onto or from a fixed
level. It is, however, possible to have a discontinuity at $x=0$.
 
In the lattice case, for $k\leq 0$,
\begin{equation}
   \label{e:HandGbis}
H(k) = G^{-k} H.
\end{equation}

\begin{prop} \label{prop:Hx}

{Assume that the process is transient.}
\begin{enumerate}
\item In the lattice case, for any $z\in {\cal D}_F$, the series 
\begin{equation}
\sum_{k\in\mathbb Z} z^{k}H(k)
\label{eq:lattHseries}
\end{equation}
 converges if and only if all the eigenvalues of $F(z)/z$ have negative real parts. If this is the case
\begin{equation}
\label{eq:lattscale}
\sum_{k\in\mathbb Z} z^{k}H(k)=-zF(z)^{-1}.
\end{equation}
\item In the non-lattice case, for any $\a\in {\cal D}_F$ the integral 
\begin{equation}
\int_{\mathbb R} e^{\a x}H(x)\D x
\label{eq:nonlattHint}
\end{equation}
converges if and only if all the eigenvalues of $F(\a)$ have negative real parts, in which case
\begin{equation}
\label{eq:nonlattscale}
\int_{\mathbb R} e^{\a x}H(x)\D x=-F(\a)^{-1}.
\end{equation}
\end{enumerate}
\end{prop}
\begin{proof}
  In the lattice case, for any $i \in E$ and $z\in {\cal D}_F$, we can take
  $f(x,j) = z^x \1{j=i}$ in~\eqref{eq:occ_density}, and then
\[\int_0^\infty z^{X_s}\1{J_s=i}\D s= \sum_{k\in\mathbb Z} z^k L_{i}(k,\infty).\]
Upon taking expectations and using  (\ref{eq:EXJ_lattice}) this yields
\[\int_0^\infty e^{sF(z)/z}\D s = \sum_{k\in\mathbb Z} z^k \e [L(k,\infty)]\]
in the obvious matrix notation. The left-hand side, and therefore the right-hand side, converges if and only if all the eigenvalues of $F(z)/z$ have negative real parts. When it converges, the integral is given by $-zF(z)^{-1}$.
In the non-lattice case this argument becomes
\[\int_{\mathbb R} e^{\a x}\e [L(x,\infty)]\D x=\int_0^\infty \e[ e^{\a X_t};J_t]\D t=\int_0^\infty e^{F(\a)t}\D t=-F(\a)^{-1}.\]
\end{proof}

\begin{remark}
\label{rem:4.2}\rm

Proposition \ref{prop:Hx} results from a nice use of the
  equivalence (\ref{eq:occ_density}) between spatial and time
  integrals of MAPs.  However, its proof depends on the hypothesis
  that the eigenvalues of $F(z)/z$ (in the lattice case) and $F(\a)$ (in
  the non-lattice case), all have negative real parts. 
 It is not always the case
that this hypothesis is satisfied for any $z$ in the domain 
$|z| < 1$ or $\a$ in  $\Re(\a)\leq 0$, respectively.  
It is for this reason that we extended the domains of $F(z)/z$ and
$F(\a)$ as defined in (\ref{eq:F_lattice}) and (\ref{eq:F}) beyond
$|z| < 1$ and $\Re(\a)\leq 0$.


Two very simple examples where it is necessary to do this are a
birth and death process with negative drift and a Brownian motion with
negative drift.
In the first example, there is no $z$ with $|z|<1$ that
  satisfies this hypothesis. But there are intervals of $z \in
(1,\infty)$ where all the eigenvalues of $F(z)/z$ have negative real
parts.  For the second example similarly, there is no $\a$ with $\Re(\a)\leq 0$
that satisfies the assumption in Proposition \ref{prop:Hx} but there is
an interval of $\a> 0$ where all the eigenvalues of 
$F(\a)$ have negative real parts. 


There are also cases where there is no value of $z$ (or $\a$) at all for which the convergence conditions of (\ref{eq:lattHseries}), respectively (\ref{eq:nonlattHint}) hold.  In both the lattice and non-lattice cases, this
occurs when the drift $\mu$ is negative  and the distribution of the jump sizes $B_i$ and $U_{ij}$ decays slower than geometrically.  

In such circumstances, we might think that the artificial killing mechanism can help. For example, we could introduce killing, say with constant rate $q$, which guarantees that there is an interval of convergence for the series (\ref{eq:lattHseries}), respectively (\ref{eq:nonlattHint}), use analytic continuation to extend the domain of (\ref{eq:lattHseries}) or (\ref{eq:nonlattHint}), and then let $q$ tend to 0. However, we have to
be careful, because this still might not increase the scope of Proposition~\ref{prop:Hx}.  In Appendix \ref{a:killing}, we present a detailed discussion of this situation, dealing with the lattice case. An analogous argument holds in the non-lattice case.

\end{remark}

Following from Proposition \ref{prop:Hx}, we can conclude a number of
things about the decay rate of $\p[J_{\tau_k}]$. The next corollary
applies to lattice processes with $\mu <0$. Analogous conclusions hold in the case where $\mu>0$ and for non-lattice processes. For $k\geq 0$. we denote the maximum eigenvalue of $\p[J_{\tau_k}]$ by $\xi_k$ and write
\[
\xi^* = \limsup_{k \to \infty} \sqrt[k]{\xi_k}.
\]
Furthermore, let $\widehat {\cal D}_F$  be  the subset of $z \in {\cal D}_F$ where all the
eigenvalues of $F(z)/z$ have negative real parts and let $\phi= \sup \{|z|: z \in \widehat {\cal D}_F\}$.

\begin{cor}
     \label{cor:decay} \ 
In the lattice case, assume that the process is non-defective. If $\mu < 0$, then
\begin{enumerate}
\item 
$\widehat {\cal D}_F $ is a subset of $\{z: |z| > 1\}$;
\item 
if $\widehat {\cal D}_F \not=\emptyset$, then $\xi^* = \phi^{-1} < 1$;
\item 
if $\widehat {\cal D}_F=\emptyset$, then $\xi^*= 1$.
\end{enumerate}
\end{cor}

\begin{proof}
By Proposition \ref{prop:Hx}, the series
$\sum_{k\in\mathbb Z} z^k H(k)$ converges
if and only if $z \in \widehat {\cal D}_F$,
and by~(\ref{eq:HxH}),
\begin{equation}
   \label{eq:splitsumb}
\sum_{k\in\mathbb Z} z^k H(k) = \sum_{k\in\mathbb Z} z^{k}
\p[J_{\tau_{k}}]H
= \sum_{k=1}^\infty z^{-k} G^k H + \sum_{k=0}^\infty z^{k} \p[J_{\tau_{k}}]H.
\end{equation}
Now, for a non-defective process with $\mu< 0$, $G$ is stochastic and
so the first sum on the right-hand side of (\ref{eq:splitsumb}) does not converge for any $z$
with $|z| \leq 1$. This proves the first claim.

The second sum on the right-hand side of (\ref{eq:splitsumb}) converges if and only if the
series $\sum_{k=0}^\infty z^{k} \xi_k$ converges. The radius of
convergence of this series is $1/\xi^*$.  If $\widehat{\cal D}_F \not= \emptyset$, then a necessary and sufficient condition for such convergence is that $|z| \leq 1/\xi^*$.  As this holds for any $z$ in $\widehat{\cal D}_F$, the second claim is proved.

The second sum in (\ref{eq:splitsumb})  converges for any $z$ in the open unit disk $\DD$ as $\p[J_{\tau_k}]$ is a stochastic or a substochastic matrix.  As a consequence, $\sum_{k=0}^\infty z \xi_k$ converges for any $z \in \DD$, which means that $\xi^* \leq 1$.  On the other hand, as $\widehat{\cal D}_F$ is empty, the second sum in (\ref{eq:splitsumb}) diverges for all $|z|>1$, from which we conclude that that $\xi^* \geq 1$ as well.  This proves the third claim.
\end{proof}


Using the fact that the matrix $H(\cdot)$ has, for negative arguments,
a simple expression in terms of $G$ and $H$, in the next
proposition we derive transforms of $\p[J_{\tau_k}]$ and $\p[J_{\tau_x}]$ over just the positive arguments.

\begin{prop}
     \label{cor:transform} \ 
{Assume that the process is transient.}
\begin{enumerate}
\item In the lattice case, 
\begin{equation}
\label{eq:latticeh}
\sum_{k=1}^\infty z^{k}\p[J_{\tau_k}]H+z(z\matI-G)^{-1}H=-zF(z)^{-1},
\end{equation}
for any $z \in \DD \setminus \sp(G)$,  where $\sp(G)$ is the set of eigenvalues of $G$.

\item In the non-lattice case, 
\begin{equation}
\label{eq:nonlatticeh}
\int_0^\infty e^{\a x}\p[J_{\tau_x}]\D xH+(\a\matI-G)^{-1}H=-F(\a)^{-1},
\end{equation}
for any $\a\in \CC_< \setminus \sp(G)$, where $\CC_<$ is the set of
complex numbers with strictly negative real part.
\end{enumerate}
The matrix $H$ is non-singular in both cases.
\end{prop}
\begin{proof}
As seen in the proof of Corollary \ref{cor:decay}, in the lattice case, the series
  $\sum_{k=1}^\infty z^{k}H(k) = \sum_{k=1}^\infty
  z^{k}\p[J_{\tau_k}]H$  converges for any $z$ in $\DD$.
Furthermore, the series
\[\sum_{k=-\infty}^0 z^kH(k)=\sum_{k=0}^\infty z^{-k}G^{k}H\]
converges to $(I - z^{-1}G)^{-1}H$  for any $z$ with $|z| > \rho(G)$, where
$\rho(G)$ is the spectral radius of $G$.

If $G$ is substochastic, $\rho(G) < 1$, and 
$\sum_{k \in \ZZ} z^k H(k)$ converges in $\{z: \rho(G) < |z| < 1\}$,
so (\ref{eq:latticeh}) follows from Proposition \ref{prop:Hx} in
that annulus.  Since $zI-G$, $F(z)$ and
$\sum_{k=1}^\infty z^{k}\p[J_{\tau_k}]$ are analytic in $\DD$, the
proof of (\ref{eq:latticeh}) in the whole of  $\DD \setminus \sp(G)$ follows by
analytic continuation.

If $G$ is stochastic, we introduce killing at the constant rate $q$
and we find by the same argument that
\begin{equation}
\label{eq:latticehq}
\sum_{k=1}^\infty z^{k}\p[J_{\tau_k};q]H_q+z(z\matI-G_q)^{-1}H_q=-zF_q(z)^{-1},
\end{equation}
in $\DD \setminus \sp(G_q)$.  By continuity of the spectrum of
$G_q$, for any given $z$, (\ref{eq:latticeh}) results from
(\ref{eq:latticehq}) in the limit as  $q \rightarrow 0$.

Finally, to see that $H$ is non-singular, we rewrite (\ref{eq:latticeh}) as
\[
-z^{-1}F(z) \left[\sum_{k=1}^\infty z^{k}\p[J_{\tau_k}]+z(z\matI-G)^{-1}\right]H=I.
\]
Similar reasoning yields the result in the non-lattice case. 
\end{proof}

\section{Returns to zero and the matrix $R$}
\label{s:R}

\subsection{Discrete-time processes: a digression} \label{sec:discrete}

In discrete-time, the occupation time $\widetilde L_j(x,n)$ counts the
number of visits to state $(x,j)$ up to step $n-1$, and we have  for the Markov
chain embedded at jump epochs (see Section \ref{sec:latt1}) 
\[
\widetilde L_{j}(x,n) =\sum_{k=0}^{n-1}\1{\widetilde X_k=x,\widetilde
  J_k=j}.
\]
Define the matrices $\widetilde R(k)$, $k \geq 1$, such that
\[
\widetilde R_{ij} (k) = \E_i[\wL_j(-k,\wtau_0) ], \qquad  i, j \in E,
\]
where $\wtau_0 =\inf\{n\geq 1:X_n=0\}$ is the first return time to level 0. Starting from level 0 in phase $i$ at time 0, $\widetilde R_{ij}(k)$ is the expected number of visits to $(-k,j)$ before the next time the process visits level 0.
It so happens that $\widetilde R(k) = \widetilde R^k$, where $\widetilde R$ is the minimal non-negative solution of the equation%

\begin{equation}
   \label{e:rtilde}
\widetilde R= \widetilde A_{-1} + \widetilde R \widetilde A_0 + \widetilde R^2 \widetilde A_1 +  \cdots
\end{equation}
The physical justification for this is that the term
$\widetilde R^k \widetilde A_{k-1}$ on the right-hand side is the
expected number of visits to level $-1$ that immediately follow a
direct jump from level $-k$.  For details, we refer to Neuts~\cite[Chapter
1]{neuts_book};  the argument there is for Markov chains which are
skip-free upward, but it is easily adapted.


\subsection{The lattice case} \label{sec:Rlattice}

In continuous time, \eqref{e:rtilde} is equivalent to%
\footnote{The notation $R$ has been extensively used in the
  Matrix-Analytic literature for processes which are skip-free upward,
  where that matrix counts visits to $(1,j)$ instead of $(-1,j)$.
  Readers should keep it mind the different meaning given to $R$ in the
  present paper.  Our presentation for continuous-time processes is
  only slightly different from that in
  Neuts~\cite[Section~1.7]{neuts_book}}
\begin{equation}
   \label{eq:R}
A_{-1} + R A_0 + R^2 A_1 + R^3 A_2 + \cdots =0
\end{equation}
with 
\begin{equation}
   \label{e:rrtilde}
R = \Delta_A \wR \Delta_A^{-1}.
\end{equation}
Again, $R$ is the 
minimal nonnegative solution to \eqref{eq:R}, \cite[Theorem 1.7.2]{neuts_book}.
The physical meaning of $R$ is slightly different from that of $\wR$.
However, for any fixed $i$, $j$ and $k$,
\begin{align*}
R_{ij}^k  & = [A_0]_{ii}\wR_{ij}^k [A_0]_{jj}^{-1}
\\
  & = [A_0]_{ii} \E_i[L_j(-k,\tau'_0)],
\end{align*}
where $\tau'_0$ is the first jump time to level 0, that is, 
\begin{equation}
\label{e:tauprime}
\tau'_0 = \inf\{t \geq \tau^*: X_t=0\},
\end{equation}
and
\[
\tau^* = \inf\{t > 0: (X_t,J_t) \not= (X_0,J_0)\}
\]
is the first jump time of the Markov chain.
Putting these together, 
\begin{equation}
   \label{e:rlattice}
R_{ij}^k = \frac{\E_i[L_j(-k,\tau'_0)]}{\E_i[L_i(0,\tau'_0)]},
\end{equation}
which we express as the expected sojourn time in $(-k,j)$ per unit of
time spent in $(0,i)$.%
\footnote{ Actually, this interpretation holds for discrete-time
  Markov chains, if we consider that a jump occurs every  unit of time.}

We observe from (\ref{eq:G}) and (\ref{eq:R}) that $G$
and $R$ are respectively the right- and the left solutions of the same
matrix equation $F(\cdot)=0$. We shall see below that this is a
property that holds for non-lattice MAPs as well.  

\begin{rem} \em
For future reference we note that that, in the transient case,
\begin{equation}
   \label{e:rgh}
GH = HR.
\end{equation}
To see this, define $(\widetilde H)_{ij}$ to be the expected number of visits
of the
embedded jump chain
to state $(0,i)$, given that the process started in $(0,j)$.  With
this, $H=\widetilde H\Delta_A^{-1}$ and, for $k \geq 0$,
\begin{equation}
   \label{e:rghk}
G^kH = \e [L(-k,\infty)]=\widetilde H\widetilde R(k)\Delta_A^{-1}=HR^k
\end{equation}
proving \eqref{e:rgh} in particular.
\end{rem}



\subsection{The non-lattice case} \label{sec:Rnonlattice}



We want to define matrix functions $R(x)$ with properties similar to those of $R^k$ in the lattice case.  To that end, we mimic (\ref{e:rlattice}) and set
\begin{equation}
   \label{e:rnonlattice}
R_{ij}(x) = \frac{\E_i[L_j(-x,\theta')]}{\E_i[L_i(0,\theta')]},
\end{equation}
for the random variable 
\begin{equation}
\label{e:thetaprime}
\theta' = \inf\{t \geq \theta^*: X_t = 0\}
\end{equation}
where 
\[
\theta^* = \inf\{t \geq 0: J_t \not= J_0\}.
\]
is the first jump time for the phase Markov chain. 

\begin{thm}
   \label{t:semigroup}
There exists a matrix $R$ such that 
\begin{equation}
\label{e:er}
R(x)=e^{Rx}. 
\end{equation}
It is a left solution of the equation $F(\cdot) = 0$.  
\end{thm}

To prove Theorem \ref{t:semigroup}, we proceed in three steps.  To
begin with, we assume that the process $(X,J)$ is transient, so that
$H_{ij}(x) = \E_i [L_j(x, \infty)]$
is finite for all $i$, $j$ and  $x$. 

\begin{lem}
   \label{t:HandR}
{\bf (Step 1)}
If $H$ is finite, then
\begin{equation}
   \label{e:HandR}
H(-x) = H R(x)
\end{equation}
for $x \geq 0$,
\end{lem}
\begin{proof}
Using a sequence of random variables with a distribution identical to that of $\theta'$ in
(\ref{e:thetaprime}), we define a sequence $\{S_k:k \geq 0\}$ of
epochs when $X_t$ returns to 0:  $S_0=0$ and $S_{k+1}= S_k +\theta'_k$
with
\[
\theta'_k = \inf\{t \geq \theta_k^*: X_{S_k+t} = 0\}
\]
where 
\[
\theta^*_k = \inf\{t \geq 0: J_{S_k+t} \not= J_{S_k}\}.
\]
As we assume that $(X, J)$ is transient, this is a terminating process
and there is an index $K^*$, finite with probability 1, after which
$X_t$ does not return to $0$, that is, $S_{K^*}< \infty$ and
$\theta'_{K^*} = \infty$. By convention, we take $S_k=\infty$ and
$J_{S_k} = \dagger$ for $k\geq K^*+1$. 

The successive phases $\{J_{S_k}\}$ form a Markov chain on the state
space $E\cup \{\dagger\}$ with absorbing state $\dagger$ and
transition matrix  
\[
P^*_{i,j} = P_i[\theta'_0< \infty, J_{S_1}=j]
\]
for $i,j \in E$. The matrix $P^*$ is substochastic but not stochastic
and has spectral radius less than one. 
We write
\[
H_{ij}(-x) = \E_i [L_j(-x, \theta'_0)]+\E_i [L_j(-x, \infty)] - \E_i [L_j(-x, \theta'_0)]
\]
with $\E_i [L_j(-x, \infty)] - \E_i [L_j(-x, \theta'_0)] = 0$ if
$\theta'_0=\infty$.  Conditioning on $J_{S_1}$, we obtain  
\begin{align}
H_{ij}(-x) & = \E_i [L_j(-x, \theta'_0)] + \sum_{u \in E} P^*_{iu} H_{uj}(-x)
\nonumber\\
& = \sum_{k\geq 0}\sum_{u \in E}\left(P^*\right)^k_{iu} \E_u[L_j(-x, \theta'_0)]
\nonumber\\
& = \sum_{u \in E}  (I-P^*)^{-1}_{iu} \E_u[L_j(-x, \theta'_0)]
\label{e:Hofx}
\end{align}
by a standard argument, noting that $\E_u[L_j(-x, \theta'_0)]\geq 0$ for all $u$ and $j$ and the exchange of summations is justified by the monotone convergence theorem. By (\ref{e:rnonlattice}), it follows that
\begin{equation}
\label{e:hijx}
H_{ij}(-x) = \sum_{u \in E}  (I-P^*)^{-1}_{iu} \E_u[L_u(0, \theta'_0)]R_{uj}(x).
\end{equation}
Furthermore, by (\ref{e:Hofx}), 
\begin{align*}
H_{iu} & = \sum_{v \in E}  (I-P^*)^{-1}_{iv} \E_v[L_u(0, \theta'_0)]
\\
& = (I-P^*)^{-1}_{iu} \E_u[L_u(0, \theta'_0)]
\end{align*}
since $\E_v[L_u(0, \theta'_0)]=0$ for $u \not= v$. Substituting into (\ref{e:hijx}), we see that
\[
H_{ij}(-x) = \sum_{u \in E} H_{iu} R_{uj}(x),
\]
which concludes the proof.
\end{proof}

Recall from Proposition  \ref{cor:transform} that $H$ is non-singular
when $(X,J)$ is transient.  Consequently, by (\ref{e:HandR}) and
(\ref{e:HandG}), for any $x,y>0$,
\begin{equation}
   \label{e:semigrp}
H R(x+y) = H(-(x+y)) = e^{G(x+y)} H = e^{Gx} H H^{-1} e^{Gy} H = H R(x) R(y)
\end{equation}
and so $R(x+y)= R(x) R(y)$. Putting this together with the fact that $R(0)=I$, the form
(\ref{e:er}) then follows by \cite[Theorem 2.9]{enna00}. 

\paragraph{Step 2}  If the process is null-recurrent, then $H$ is not finite and we introduce an artificial killing rate $q$;  the
semi-group property follows by continuity of $R_q(x)$ as $q$ goes to 0.

\paragraph{Step 3: Time reversal}  
To prove the second statement in Theorem \ref{t:semigroup}, we proceed
with time reversal, a tool much used in the analysis of MAPs (see~\cite{asra90,rama90}).

Consider at first  a
non-defective MAP $(X,J)$ and let $\bs\pi$ be the stationary
distribution of~$J$. For an arbitrary $t>0$ and for $s\in[0,t)$, we
define a time-reversed process by
\begin{equation}
   \label{eq:timereverse}
\widehat J_s=J_{(t-s)-}\   \mbox{ and }\  \widehat X_s=X_t-X_{(t-s)-}.
\end{equation}
Assume that $J_0$ has distribution $\bs \pi$. Then $\widehat J_s$ is a stationary time-homogeneous Markov chain and $(\widehat X_s,\widehat J_s)$ is a MAP, which
is characterized by 
\begin{align}
    \label{eq:timerev}
\widehat F(\a)=\Delta_\pi^{-1}F(\a)^T\Delta_\pi
\end{align}
independently of the time horizon~$t$. 
Notice that over intervals of time where $\whJ$ is constant, and equal to $i$, say, the process $\whX$ behaves as the original $X^{(i)}$, and so $\whX$ is one-sided, like $X$.
Moreover, the stationary drift $\hat\mu$ is equal to $\mu$.

If the process is defective, $\bs\pi$ should be the stationary
distribution of the non-defective rate matrix $Q+\Delta_{\bs q}$,
where $\bs q$ is the vector of killing rates; this clearly preserves
the identity~\eqref{eq:timerev}.  Note that by reversing time twice we
obtain the original process and, furthermore, we have some freedom in
defining the reversed process, see~\cite{asm_iva}.

\begin{lem}
     \label{thm:GRH}
Let $\widehat G, \widehat H, \widehat R$ be the analogues of the
matrices $G,H,R$ for the  
process characterized by $\widehat F(\a)$ in~\eqref{eq:timerev}.  Then
\begin{gather}
\widehat G=\Dpi^{-1}R^T\Dpi,\qquad \widehat R=\Dpi^{-1}G^T\Dpi.
\label{eq:3id}\end{gather}
If the process is transient, then 
\begin{equation}\label{eq:hatH}
\widehat H=\Dpi^{-1}H^T\Dpi.
\end{equation}
The matrix $R$ is the unique left solution of $F(\cdot)=0$ in the class of matrices specified in Section~\ref{sec:nonlatt1}.
\end{lem}
\begin{proof}
  As usual, we start with the transient case and we assume that
  \eqref{eq:nonlattscale} holds for some $\alpha$ in ${\cal D}_F$.  It
  follows that for all $x\in\R$ we have 
\begin{equation}
   \label{e:HandhatH}
\widehat H(x)=\Dpi^{-1}H(x)^T\Dpi,
\end{equation}
because \eqref{eq:nonlattscale}
identifies~$H(x)$; we obtain \eqref{eq:hatH} by setting $x=0$.   Combining this identity with \eqref{e:HandR}, we
readily obtain
\begin{align*}
H e^{Rx} & = \Delta_\pi^{-1}  \whH(-x)^T \Delta_\pi && \mbox{by
                                                          (\ref{e:HandhatH})}
\\
 & = \Delta_\pi^{-1}  (e^{\whG x} \whH)^T \Delta_\pi && \mbox{by
                                                          (\ref{e:HandG})}
\\
 & = \Delta_\pi^{-1}  \whH^T (e^{\whG x})^T \Delta_\pi
\\
 & = H \Delta_\pi^{-1}  (e^{\whG x})^T \Delta_\pi
\end{align*}
for any $x\geq 0$.  As $H$ is nonsingular, this proves that $\widehat
G=\Dpi^{-1}R^T\Dpi$ and one shows by a similar argument that $\widehat
R=\Dpi^{-1}G^T\Dpi$. 

The non-defective case with $\mu = 0$ can be treated by introducing an
artificial killing rate $q$, concluding that \eqref{eq:3id} holds for the processes
with killing, and  letting the killing rates go to~0.

Finally, we recall from Section~\ref{s:ladder} that $\widehat G$ is
the unique right solution of $\widehat F(\cdot)=0$ in a certain
class. Using the representations of these quantities in terms of the
original process, one readily verifies that $R$ is the left solution of $F(\cdot)=0$, unique in the given class. 
\end{proof}

\begin{remark}\label{rem:S1} \rm
We can replace the stopping time $\theta^\prime$ in the definition (\ref{e:rnonlattice}) of $R(x)$ with any other stopping time 
 \[\eta\leq \overline\eta_0=\inf\{t \geq 0: X_t = 0, J_t \not= J_0\}\] that is such that
$X_{\eta}=0$ when $\eta<\infty$, and $\e_i L_{i}(0,\eta)>0$ for all~$i$. It is clear that any such stopping time leads to the same matrix~$R(x)$.
For instance, when the processes $X^{(i)}$ have bounded variation for
all $i$, with linear drifts $d_i<0$, we may choose
$\eta=\inf\{t>0: X(t)=0\}$ as the first return to~0.  
It follows
from~\eqref{eq:Ldefbounded} that
\[R_{ij}(x)=\frac{|d_i|}{|d_j|}\e_i\#\{t\in[0,\eta):X_t=-x,J_t=j\},\]
which closely resembles the discrete-time lattice case.
\end{remark}

\begin{rem} \em
   \label{r:GHR}
If $H$ is finite, then by \eqref{e:HandG}, \eqref{e:er} and
\eqref{e:HandR},
\begin{equation}
   \label{e:GHR}
G H = H R,
\end{equation}
identical to \eqref{e:rgh} and we conclude, using \eqref{eq:3id} that
the four matrices  $G$, $R$, $\whG$ and $\whR$ are algebraically similar.
\end{rem}


\begin{rem} \em
It is interesting to consider whether it is possible to prove Theorem \ref{t:semigroup} without writing $R$ in terms of $G$ and using the analogous properties of $G$. 

The semi-group property
\[
R_{ij}(x+y)= \sum_{u \in E} R_{iu}(x) R_{uj}(y)
\]
is equivalent to 
\[
\E_i [L_j(-x-y, \theta'] = \sum_{u \in \E} \E_i [L_u(-x, \theta']  R_{uj}(y),
\]
which indicates that the local time spent in $(-x-y,j)$ is the sum of the local time spent in the various phases at level $-x$ multiplied by the local time spent $y$ units lower per unit of local time spent at level $-x$. 

Since the process is skip-free downward, starting from level 0, it is obviously necessary for the process to move through level $-x$ before reaching level $-x-y$. Therefore, it is tempting to repeat the argument used in~\cite[Theorem 6.2.7]{lara99} for lattice processes, and  decompose the interval $(0,\theta')$ into subintervals of visits to level $-x$.  Unfortunately, a formal proof of the semi-group property starting from that observation has eluded us so far.
\end{rem}




\section{Two-sided exit and matrix scale functions }
\label{s:two-sided}

Two-sided exit theory is concerned with a study of events
$[\tau_{-a}<\tau_b^\u]$ for $a,b\geq 0$, where 
\begin{equation}
   \label{eq:tauxplus}
\tau_x^+=\inf\{t\geq 0: X_t\geq x\}
\end{equation}
is the first passage time over a non-negative level $x$.  For the
non-lattice case, so-called matrix-valued scale functions $W(x)$ play
an important role, and we recall relevant results below.  Later in
this section, we derive their counterparts for the lattice case.

\subsection{Non-lattice MAPs.} \label{sec:nonlatticescale}

The main result in \cite{ivanovs_palmowski} states that, for $a,b\geq 0$ with $a+b>0$, there exists a continuous matrix-valued function $W(x), x\geq 0$, nonsingular for $x>0$, such that 
\begin{equation}
   \label{e:latDab}
\p[\tau_{-a}<\tau_b^+,J_{\tau_{-a}}]=W(b)W(a+b)^{-1}.
\end{equation}
The function $W(x)$ is such that 
\begin{equation}
   \label{e:W}
\int_0^\infty e^{\a x}W(x)\D x=F(\a)^{-1}
\end{equation}
for $\a< \min\{\Re(\lambda):\lambda \in \sp(G) \}$.
Furthermore,
\begin{equation}
   \label{eq:WTheta}
W(x)=e^{-G x}\Theta(x),
\end{equation}
where
\begin{equation}
   \label{e:Theta}
\Theta(x)=\e [L(0,\tau_{-x}) ] \quad \mbox{for $x>0$,} 
\end{equation}
and $\Theta(0)=\lim_{x\downarrow 0}\Theta(x)$.

This identity implies that $W(0)=\Theta(0)$ and $\Theta_{ij}(0)=0$ unless $i=j$ and $X_i$ is a bounded variation process, in which case $\Theta_{ii}(0)=1/|d_i|$, see (\ref{eq:Ldefbounded}). This observation is consistent with the fact that $\p_i[\tau_0^+=0]$ is either~0 or~1 according to $X^{(i)}$ having bounded or unbounded variation, see~\cite[Thm.\ 6.2]{kyprianou}.

If the process is transient then, by decomposing the local time at level 0 into segments that occur before and after $X_t$ first visits level $-x$,  we see that
\begin{equation}
   \label{e:thetaH}
\Theta(x)=H-e^{Gx}\p[J_{\tau_x}]H
\end{equation}
 and, for $x > 0$, 
\begin{equation}
   \label{eq:Jtaux}
\p[J_{\tau_x}]=e^{-G x}-W(x)H^{-1},
\end{equation}
providing a simple expression for the phase at the first hitting time of a positive level $x>0$. The identity (\ref{eq:Jtaux}) can be generalized to all $x\in\mathbb R$ by putting $W(x)=0$ for $x<0$ and redefining $\tau_0$ to be the first hitting time of~$0$.

\subsection{Lattice MAPs, arbitrary $A_{-1}$} \label{sec:latticescalegenA-1}

In the lattice case, the matrix $A_{-1}$ of transition rates from level $0$ to level $-1$ plays a crucial role. It is tempting to assume that this matrix is nonsingular, in which case an analysis that mirrors that of the scalar case discussed in \cite{avvi19} is essentially feasible. However, this is a restrictive assumption from an applications point of view. For example, in the MAP corresponding to an $M/Erlang_k/1$ queue, the only transitions that decrease the level and have a positive rate move from phase $k$ to phase 1. More generally, $A_{-1}$ for an $M/Ph/1$ queue always has rank 1.

In this section, we start with lattice MAPs in which we make no assumption about the nonsingularity of $A_{-1}$. We pursue a dual objective: on the one hand we look for an expression structurally similar to (\ref{e:latDab}), and on the other hand we look for a connection with matrix functions that satisfy an equation similar to (\ref{e:W}). The first objective leads
us to Lemma~\ref{t:lattA}, while the second leads us to Lemmas \ref{t:lattB} and \ref{t:recurscale} which do require the assumption that $A_{-1}$ is nonsingular.  

For $a, b \geq 0$, define
\begin{equation}
\label{e:Dab}
D_{a,b} = \P[\tau_{-a} < \tau_b^\u, J_{\tau_{-a}} ], 
\end{equation}
and, for $m \geq 1$, define the matrix
\[
\Xi(m)= \E[L(0, \tau_m^\u) ],
\]
of expected sojourn times at level 0 until the positive level $m$ is crossed.

\begin{lem}
   \label{t:lattA}

\begin{enumerate}
\item For $a \geq 0$, $b \geq 1$, 
\begin{equation}
   \label{e:dab}
D_{a,b} = \Xi(b) \hR^a (\Xi(a+b))^{-1}.
\end{equation}
\item Also, for $a \geq 1$, $D_{a,0} = 0$  and $D_{0,0} = I$.
\end{enumerate}
\end{lem}
\begin{proof}
For $a$, $b \geq 0$, define the matrices 
\[
\Xi(a, b) = \E[L(-a, \tau_b^\u)] = \E[L(-a, \tau_b) ].
\]
of expected sojourn times at level $-a$ until the higher level $b$ is
crossed. The second equality holds because the process is forced to
visit level $b$ before $-a$ if it jumps above $b$, but the first
equality serves as a reminder that the process is allowed to jump
above $b$. 

For $m \geq 1$, we have $\Xi(m)=\Xi(0,m)$ and, obviously, 
$\Xi(a,b)= D_{a,b} \Xi(a+b)$.  Also, 
\[
\Xi(a,b) = \Xi(b) \hR^a,
\]
as the first factor on the right is 
the expected time spent in level 0 before the first visit to $b$, and
the second is the expected time spent in level $-a$ per unit of time
spent at level 0.  Bringing together the two expressions for $\Xi(a,b)$, we have  
\[
D_{a,b} \Xi(a+b) = \Xi(b) \hR^a.
\]
To conclude the proof of (\ref{e:dab}), we need to show that $\Xi(m)$ is non-singular
for $m \geq 1$. 
Define $M(m)$ to be the transition probability matrix
from level 0 back to level~0 under taboo of level $m$, with entries  
\[
M(m)_{ij} = \P_i[\tau'_0 < \tau_m^\u, J_{\tau'_0} = j],
\]
where $\tau'_0$ is the first jump time to level 0. Since the process is irreducible, $M(m)$ is a substochastic matrix that is not
stochastic, and its spectral radius is strictly less than one. Thus its power series is convergent and
\[
\Xi(m)= \sum_{\nu \geq 0} M(m)^\nu \Delta_A^{-1}= (I - M(m))^{-1}\Delta_A^{-1},
\]
which concludes the proof.
\end{proof}

Lemma \ref{t:xi} below gives an expression for the matrices $\Xi(m)$.
\begin{lem}
   \label{t:xi}
For $m \geq 1$, the expected sojourn time in level 0 under taboo of level $m$ is given by
\begin{equation}
   \label{e:xi}
\Xi(m) = \sum_{0 \leq \nu \leq m-1} \Phi(\nu) \hR^\nu,
\end{equation}
where 
\[
\Phi(m) = \E[L(m, \tau_{m+1}^\u)]
\]
is the expected sojourn time in level $m$, starting from level 0, under taboo of level $m+1$.
\end{lem}

\begin{proof}
As $R^{m-1}$ is equal to the matrix of expected sojourn times in level 0 per
unit of time spent in level $m-1$,  it is easy to see that, for $m \geq 1$, 
\[
  \E[L(0, \tau_m^\u) ] = 
  \E[L(0, \tau_{m-1}^\u) ] 
  + \E[L(m-1, \tau_m^\u) ]
  R^{m-1},
\]
that is,
\begin{align}
   \label{e:lotau}
\Xi(m) = \Xi(m-1) + \Phi(m-1) R^{m-1}.
\end{align}
By successive substitution, we find that
\[
\Xi(m)=\Xi(1)+ \Phi(1) \hR + \cdots + \Phi(m-1) \hR^{m-1},
\]
which proves the claim since $\Xi(1)=\Phi(0)$ by definition.
\end{proof}

To provide expressions for $\Phi(m)$, we use the transition matrices
\begin{equation}
\label{e:Aistar}
A_i^\u = \P[\tau_0^\u = \tau_i, J_{\tau_0^\u}]
\end{equation}
of taboo probabilities that, starting from level 0, the first nonnegative level is level $i$. Equivalently, these are the overshoot probabilities that the process moves to level $i$ when it first visits the nonnegative levels.  We have defined in (\ref{eq:jumpchain}) the transition matrices $\wA_i$ of the jump chain and a simple argument proves that, for $i \geq 0$,
\begin{equation}
   \label{e:astar}
A_i^\u = \sum_{\nu \geq 0}  \wR^\nu \wA_{\nu+i}
\qquad \mbox{or}  \qquad 
A_i^\u = \wA_i + \wR A_{i+1}^\u
\end{equation}
where $\wR$, the minimal nonnegative solution to equation (\ref{e:rtilde}), is the matrix of expected numbers of visits to level $-1$ under taboo of level 0.  

\begin{lem}
   \label{t:vm}
The expected sojourn times $\Phi(m)$ in level $m$ under taboo of level $m+1$ are given by
\begin{align}
   \label{e:xip0}
\Phi(0) & = (I- A_0^\u)^{-1}  \Delta_A^{-1}
\\
  \nonumber
\Phi(n) & = (I- A_0^\u)^{-1}  \sum_{1 \leq \nu \leq n}  A_\nu^\u \Phi(n-\nu),
\end{align}
for $n \geq 1$.
\end{lem}

\begin{proof}
We condition on the first jump to a nonnegative level and write
\begin{align}
   \label{e:vo}
\Phi(0) & = \Delta_A^{-1} + A_0^\u \Phi(0) \\
   \nonumber
\Phi(n) & = A_0^\u \Phi(n) + A_1^\u \Phi(n-1)+ \cdots + A_n^\u \Phi(0),
\end{align}
from which the lemma follows --- the matrix $I-A_0^\u$ is non singular because the process is irreducible.
\end{proof}

The next lemma gives a simple relation between the matrices $\Xi(k)$
and $\Theta(k)$. It holds for null recurrent as well as transient processes.
\begin{lem}
   \label{t:thetaxi}
For all $k \geq 1$
\begin{equation}
   \label{e:thetaxi}
\Theta(k) R^k = G^k \Xi(k).
\end{equation}
\end{lem}
\begin{proof}
Define $\delta = \inf\{t > \tau_{-k} : X_t = 0\}$, that is, assuming
the process starts in level 0, $\delta$ is the first return time to
level 0 after the first visit to level $-k$, with $\delta = \infty$ if
no such visit exists.  We express the expected time in
level $-k$ during the interval $(0, \delta)$ in two ways.

On the one hand,
\[
\E[L(-k, \delta)] = G^k \E[L(0,\tau^{\u}_k)]
\]
by the Markov property.  On the other hand, 
\[
\E[L(-k, \delta)] = \E[L(0, \delta)]  R^k
\]
since $R^k$ is the expected time in level $-k$ per unit of time in
level 0, and $\E[L(0, \delta)]  = \E[L(0, \tau_{-k})]$ since the
process does not spend any time in level 0 during the interval
$(\tau_{-k}, \delta)$.  Thus,
\[
\E[L(0, \tau_{-k})] R^k = G^k \E[L(0,\tau^{\u}_k)],
\]
which concludes the proof.
\end{proof}

\begin{rem} \em
If the process $(X,J)$ is transient, then
\begin{equation}
   \label{e:xiG}
\Xi(m) = H -\P[J_{\tau_m}] G^m H;
\end{equation}
indeed, the second term in the right-hand side is equal to the
expected time in level 0 after 
the first visit to level $m$.  This relation is similar to
(\ref{e:thetaH}).  By (\ref{e:rghk}), we can also write 
\begin{equation}
   \label{e:xiH}
\Xi(m) = H -\P[J_{\tau_m}] H R^m.
\end{equation}
\end{rem}

\begin{rem} \em In principle, the functions $\Xi(m)$ are well-determined: by (\ref{e:astar}), the matrices $A^\u_i$ are the
  solution of the system
\[
\vligne{A^\u_0 \\ A^\u_1 \\ A^\u_2 \\ \vdots}
=
\vligne{I  & -\wR \\ & I  & -\wR \\ & & I  & \ddots \\ & & & \ddots }
\vligne{\wA_0 \\ \wA_1 \\ \wA_2 \\ \vdots}
\]
and once this system is solved, the $\Phi(m)$s and the $\Xi(m)$s can be computed recursively by (\ref{e:xip0}) and (\ref{e:xi}). 
\end{rem}


\subsection{Lattice MAPs, nonsingular
  $A_{-1}$.} 
         \label{sec:latticescalenonsingA-1}   

If the matrix $A_{-1}$ is nonsingular, we can define a matrix function
$W(k)$, $k \geq 0$, {which satisfies an equation similar to
  (\ref{e:W}), and which leads to an
expression for  $D_{a,b}$  identical to \eqref{e:latDab}.
This function, clearly, is the analogue for lattice MAPs of
the scale matrix function of Section~\ref{sec:nonlatticescale}.
}

\begin{lem}
   \label{t:lattB}
If $A_{-1}$ is nonsingular, then
\begin{equation}
   \label{e:dabalt}
D_{a,b} =    (\Xi(b) R^{-b}) (\Xi(a+b) R^{-(a+b)})^{-1}
\end{equation}
for $a$, $b \geq 0$.
\end{lem}
\begin{proof}
It is well-known that $R = A_{-1} S$, where $S_{ij}$ is the expected
sojourn time in $(-1,j)$, starting from $(-1,i)$, before the first
visit to a non-negative level (see \cite[Theorem 3.1]{latou93}) and $S$ is
nonsingular.
Therefore, $R$  is nonsingular if and only if $A_{-1}$ is nonsingular
and (\ref{e:dabalt}) is a direct consequence of (\ref{e:dab}).  
\end{proof}

\begin{lem}
   \label{t:recurscale}
If $A_{-1}$ is nonsingular, the matrices $W(k) = \Xi(k) R^{-k}$ satisfy the recursion
\begin{align}
  \nonumber
W(1) & = (A_{-1})^{-1}  \\
   \label{e:wk}
W(k+1) & = - (A_{-1})^{-1}\sum_{1 \leq \nu \leq k} A_{k-\nu} W(\nu).
\end{align}
Furthermore,
\begin{equation}
\label{eq:discreteinvgen}
\sum_{k\geq 1} W(k) z^k = z F(z)^{-1}
\end{equation}
for $|z| < \min\{|\lambda|: \lambda \in \sp(G)\}$.
\end{lem}

\begin{proof}
We observed above that (\ref{e:xi}) implies that $\Xi(1)= \Phi(0)$, and so
\[
W(1) = \Phi(0) R^{-1}.
\]
Then, by (\ref{e:astar}) and (\ref{e:xip0}),
\begin{align*}
W(1) & = (I-\sum_{\nu \geq 0} \wR^\nu \wA_\nu)^{-1} \Delta_A^{-1} R^{-1} \\
& = (A_{-1})^{-1},
\end{align*}
where the final equation follows from (\ref{eq:jumpchain}) and (\ref{e:rrtilde}).

Also, by conditioning on the first transition of the jump chain, we
see that
\[
D_{1,k} = \wA_{-1} + \wA_0 D_{1,k} + \wA_1 D_{2,k-1} + \ddots +
\wA_{k-1}  D_{k,1}.
\]
Using (\ref{e:dabalt}), this is equivalent to $W(k) = \sum_{-1 \leq
  \nu \leq k-1} \wA_\nu W(k-\nu)$, from which (\ref{e:wk}) readily
follows.

By definition,
$\Phi(m)= \E[L(m,\tau_{m+1}^\u)]$, and $L(m,\tau_{m+1}^\u) > 0$ only
if the process reaches level $m$ before any higher level.  If this
happens, one may re-label the levels and count the expected time
spent in the new level 0 before the new level~1.  Thus, for $m \geq 1$,
\begin{equation}
   \label{e:phim}
\Phi(m) = \P[\tau_m = \tau_m^\u, J_{\tau_m}]
\Phi(0).
\end{equation}
  The series $\Upsilon(z) = \sum_{n
  \geq 0} z^n \Phi(n)$ converges for $z$ in the open unit disc $\DD$ since the matrices
$\P[\tau_n = \tau_n^\u, J_{\tau_n}]$ are bounded.
 Using Lemma \ref{t:vm}, we have
\[
(I-A_0^\u) \Upsilon(z) = \Delta_A^{-1} + \sum_{n \geq 1} \sum_{1 \leq \nu \leq n} A_\nu^\u \Phi(n-\nu) z^n,
\]
or
\begin{align*}
\Upsilon(z) & = \Delta_A^{-1} + A_0^\u \Upsilon(z) + \sum_{\nu \geq 1} A_\nu^\u z^\nu \sum_{n \geq
       \nu} \Phi(n-\nu) z^{n-\nu} \\
  & = \Delta_A^{-1} + A^\u(z) \Upsilon(z),
\end{align*}
where $A^\u(z)=\sum_{\nu \geq 0} A_\nu^\u z^\nu$ also converges in $\DD$.  Thus, 
\begin{equation}
   \label{e:av}
(I-A^\u(z)) \Upsilon(z) = \Delta_A^{-1}.
\end{equation}

By an argument  similar to \cite[Theorem 5.8]{blm05}, we show
that $\wF(z)-z I=(\wR - z I) (I - A^\u(z))$ and, together with
(\ref{e:av}),  this proves that
\[
\Upsilon(z) =  (\wF(z) - zI)^{-1} (\wR -zI) \Delta_A^{-1}
\]
if $\wF(z) - zI$ is nonsingular, that is, if $z$ is not an eigenvalue
of $G$, by \cite[Equation 4.22]{blm05}.

We replace the $\wA_i$s by their expression (\ref{eq:jumpchain}) and $\wR$ by its expression from
(\ref{e:rrtilde}) and obtain
\begin{equation}
   \label{e:thetafz}
\Upsilon(z) = F(z)^{-1} (R-zI)
\end{equation}
if $z \not\in \sp(G)$.
%
To prove (\ref{eq:discreteinvgen}), we
successively write
\begin{align*}
\sum_{k \geq 1} W(k) z^k  & = \sum_{k \geq 1} \Xi(k) R^{-k}z^k  \\
  &  = \sum_{k \geq 1} \sum_{0 \leq \nu \leq k-1} \Phi(\nu)
    R^{\nu-k}z^k    \qquad \mbox{by \eqref{e:xi}}\\
  &  = \sum_{\nu \geq 0}  \Phi(\nu) z^\nu  \sum_{k \geq 1} R^{-k} z^k
  \\
  &  = \Upsilon(z) ((I-z R^{-1})^{-1} -I) 
\end{align*}
if $|z \rho(R^{-1})| < 1$, that is, if $|z|$ is less than the smallest
absolute value of the eigenvalues of $R$.  The eigenvalues of $R$
coincide with those of $G$ and, using \eqref{e:thetafz},  we
easily complete the proof of (\ref{eq:discreteinvgen}).
\end{proof}

It is clear from (\ref{e:wk}) that the matrices $W(k)$ are of mixed
sign, and one might run into numerical instabilities if one were to
use the recursive scheme for numerical purposes.

If $A_{-1}$ is nonsingular, then we may use the expected time at level
0 under taboo of {\em lower levels}, and establish a direct connection
to the expressions (\ref{e:latDab}, \ref{eq:WTheta}) for non-lattice processes.

\begin{lem} 
   \label{t:lattC}
If $A_{-1}$ is nonsingular, then
\begin{equation}
   \label{e:dabter}
D_{a,b} = (G^{-b} \Theta(b))\, (G^{-a-b} \Theta(a+b))^{-1},  \qquad
\mbox{for $a \geq 0$, $b \geq 1$,}
\end{equation}
where $\Theta(k)$ is defined in (\ref{e:Theta}).
\end{lem}
\begin{proof}
As $A_{-1}$ is nonsingular, both $G$ and $R$ are nonsingular and it follows from Lemma~\ref{t:thetaxi} that
$G^{-k} \Theta(k) = \Xi(k) R^{-k}$.  With this,  (\ref{e:dabter}) readily
follows from (\ref{e:dabalt}).
\end{proof}






\begin{remark} \em
\label{rem:RV}
Observe that $\sum_{m=0}^\infty \Phi(m)$ is the expected time the process spends at its running maximum. By (\ref{e:thetafz}), for a defective process this quantity should be finite, and equal to $A^{-1}(R-\matI)$.  In the non-lattice case the corresponding quantity is given by the potential density at $0$ of the process reflected at~0 from above, which is $Q^{-1}R$ see~\cite[Equation (24)]{ivanovs_potential}. These identities provide yet another possible interpretation of~$R$.
\end{remark}

\section{Applications}
   \label{s:applications}
\subsection{Creeping} \label{sec:creeping}

A \levy process $X$ is said to {\it creep} over a level $x>0$
if 
\[\p[X_{\tau^+_x}=x , J_{\tau_x}]=\p[\tau_x=\tau_x^+, J_{\tau_x}]>0,\] 
see
~\cite[{\S}7.5]{kyprianou}. This means that $X$ can upcross $x$
continuously.  A related problem would be to consider
$\p[\tau_x<\tau_{x+y}^+ , J_{\tau_x}]$ for $x,y>0$, that is,  the probability of
hitting a level $x$ before crossing over $x+y$. 

\paragraph{The lattice case.}

For $m \geq 1$, we repeat the argument in Lemma \ref{t:vm}, conditioning on the first
jump to a nonnegative level, and obtain
\begin{align*}
\p[\tau_0=\tau_0^+, J_{\tau_0}]  & = A_0^\u  \\
\p[\tau_m=\tau_m^+, J_{\tau_m}]  & =  (I-A_0^\u )^{-1} \sum_{1 \leq
                                   \nu \leq m} A_\nu^\u \p[\tau_{m-\nu}=\tau_{m-\nu}^+, J_{\tau_{m-\nu}}].
\end{align*}
By (\ref{e:xip0}), $\Phi(0)$ is nonsingular and then by (\ref{e:phim}),  for $m \geq 1$, we have
\begin{equation}
   \label{e:pv}
\P[\tau_m = \tau_m^\u, J_{\tau_m}] = \Phi(m) \Phi(0)^{-1}.
\end{equation}
By an argument similar to (\ref{e:phim}),  the expected time at level $m$ before $\tau^\u_{m+l}$ may
be expressed as
\begin{equation}
   \label{e:lim}
\E[L(m, \tau^\u_{m+l})] =
\P[\tau_m < \tau^\u_{m+l}, J_{\tau_m}] \Xi(l) 
\end{equation}
and by an argument similar to (\ref{e:xi}), as
\begin{equation}
   \label{e:limb}
\E[L(m, \tau^\u_{m+l})] =
\sum_{0 \leq i \leq
  l-1} \Phi(m+i) R^i.
\end{equation}
Thus,
\begin{equation}
   \label{e:tautaup}
\P[\tau_m < \tau^\u_{m+l}, J_{\tau_m}] = (\sum_{0 \leq i \leq
  l-1} \Phi(m+i) R^i)
\, (\sum_{0 \leq i \leq  l-1} \Phi(i) R^i)^{-1}.
\end{equation}
Finally, we may write
\[
\P[ J_{\tau_m}] = \lim_{l \rightarrow \infty} \big\{(\sum_{0 \leq i \leq
  l-1} \Phi(m+i) R^i)
\, (\sum_{0 \leq i \leq  l-1} \Phi(i) R^i)^{-1} \big\}.
\]
The limit exists since the left-hand side of (\ref{e:tautaup}) is an
increasing function of $l$ and is bounded above by 1.  If the process
is transient, then both series in the equation  above are finite, and
this becomes
\[
\P[ J_{\tau_m}]  = \E[L(m,\infty)]  \, \E[L(0,\infty)]^{-1}, 
\]
identical to (\ref{eq:HxH}).

If $A_{-1}$, and hence $R$, is nonsingular, we can write (\ref{e:tautaup}) in terms of the scale matrices $W(k) = \Xi(k)R^{-k}$,
\begin{align}
   \nonumber
\P[\tau_m < \tau^\u_{m+l}, J_{\tau_m}]  & = 
  (\Xi(m+l) - \Xi(m)) R^{-m}  (\Xi(l))^{-1} && \mbox{by(\ref{e:xi},
                                               \ref{e:lim}, and \ref{e:limb}),}
\\
   \label{e:PW}
  & = (W(m+l) - W(m) R^{-l}) (W(l))^{-1}.
\end{align}
In particular, 
\[
\P[\tau_m = \tau_m^\u, J_{\tau_m}] =
\P[\tau_m < \tau^\u_{m+1}, J_{\tau_m}] =
 (W(m+1) - W(m) R^{-1}) (W(1))^{-1},
\]
giving us another expression for the creeping probability.

\paragraph{The non-lattice case.}


In the non-lattice case, we start from (\ref{e:lim}) written here as
\[
\p[\tau_x<\tau_{x+y}^+,J_{\tau_x}]\Xi(y)=\e [L(x,\tau_{x+y}^+)].
\]
If the process is transient, we have 
\[
\e[L(x,\tau_{x+y}^+)]=\p[J_{\tau_x}]H-\p[J_{\tau_{x+y}}]e^{G y}H
\]
from the same argument as (\ref{e:xiH}).  Now, in the non-lattice
case, $W(x)$ is defined in (\ref{eq:WTheta}) as $W(x) = e^{-Gx}
\Theta(x)$, which is equal to $\Xi(x) e^{-Rx}$ by the non-lattice
version of (\ref{e:thetaxi}).  We use (\ref{eq:Jtaux}, \ref{e:GHR})
and eventually obtain
\begin{equation}
    \label{eq:creeping_gen}
\p[\tau_x<\tau_{x+y}^+,J_{\tau_x}]=(W(x+y)-W(x)e^{-Ry})W(y)^{-1},
\qquad x\geq 0,y>0,
\end{equation}
the same equation  as (\ref{e:PW})  in the
lattice case when $A_{-1}$ is nonsingular.

The zero-drift case can be treated by the standard limiting
argument. 


Note
that $\e[L(x,\tau_{x+y}^+)]$ is just the potential density at level $x$
for the process killed upon crossing $x+y$; its formula was identified
in~\cite[Thm.\ 1]{ivanovs_potential}. 

We conclude by providing the creeping probability. Rewrite the right hand side of~\eqref{eq:creeping_gen} as
\[
( W(x+y)-W(x)+W(x)(\matI-e^{-Ry}) ) W(y)^{-1}
\]
and let $y\downarrow 0$ to obtain
\[
\p[\tau_x=\tau_{x}^+,J_{\tau_x}]=(W'(x)+W(x)R)\lim_{y\downarrow
  0}yW(y)^{-1}, \qquad x>0,
\]
where $W'$ denotes the right derivative of~$W$, proved to exist
in~\cite{ivanovs_palmowski}.  This formula is a generalization of that in~\cite[Exercise
8.6]{kyprianou} for a spectrally-positive \levy process. If $X^{(i)}$ is of bounded variation then, under Assumption \ref{a:nosub}, $X$ cannot
creep over any level in phase~$i$, see~\cite[Thm.\
7.11]{kyprianou}. If all of the $X^{(i)}$ are of unbounded variation
then the limit reads as $W'(0)^{-1}=\frac{1}{2}\Delta_{\bs\sigma}^2$
which can be established as in the \levy case. 

\subsection{Occupation times}\label{sec: occupation}

\paragraph{The lattice case.}

\begin{lem}
   \label{t:linf}
Assume the process is transient. The expected total sojourn time in
level 0 is given by
\begin{equation}
   \label{e:l0}
\E[L(0, \infty)]  = (I - U^\u)^{-1} \Delta_A^{-1},
\end{equation}
where
\begin{equation}
   \label{e:ustar}
U^\u = \sum_{i \geq 0} A^\u_i G^i
\end{equation}
is the transition matrix from a state in level 0 back to a state in
level 0.  For $m > 0$,
\begin{align}
   \label{e:lmm}
\E[L(-m, \infty)]  & = G^m \, \E[L(0, \infty)]
\\
   \label{e:lm}
\E[L(m, \infty)]  & = (I-A^\u_0)^{-1} \big\{\sum_{1 \leq \nu \leq m} A^\u_\nu \E[L(m-\nu,
                    \infty)] 
  \\  \nonumber
 & \qquad + \sum_{\nu \geq m+1} A^\u_\nu G^{\nu-m}  \E[L(0, \infty)]\big\},
\end{align}
where the transition matrices $A^\u_i$, for $i \geq 0$, are defined in
(\ref{e:Aistar}).
\end{lem}

\begin{proof}
One readily sees from its definition (\ref{e:ustar}) that  $U^\u_{ij}$ is equal to the probability $\P_i[\tau'_0 < \infty, J_{\tau'_0}=j]$ with $\tau'_0$ defined as in (\ref{e:tauprime}). As the process is transient, $U^\u$ is substochastic but not stochastic, the inverse of $I-U^\u$ is well-defined and its entries are the expected total number of visits to level 0.  This justifies (\ref{e:l0}). 

Equation (\ref{e:lmm}) is merely a repeat of (\ref{e:HandGbis}).
For positive levels, we condition on the first visit to a nonnegative
level and write
\[
\E[L(m, \infty)] = \sum_{\nu \geq 0} A^\u_\nu \E[L(m-\nu, \infty)],
\]
from which (\ref{e:lm}) follows after simple manipulation.
\end{proof}

Several expressions are available for sojourn times in a given level
under taboo of some other levels. We only give two sets here.  The proof of the next lemma is
immediate and is omitted. 

\begin{lem}
   \label{t:nonnull}
Assume that the process is transient.  We have
\begin{align}
   \nonumber
\E[L(k,\tau^\u_m)] & = \E[L(k,\infty)] - \E[L(m, \infty)] R^{m-k}, &&
                                           k < m,
\\
   \nonumber
\E[L(k, \tau_{-l}]  & = \E[L(k, \infty)] - G^l \E[L(k+l, \infty)],  &&
                                           k > -l,
\\
\E[L(k, \tau_{-l} \wedge  \tau^\u_m)]  & =
   \E[L(k,\tau^\u_m)] - D_{l,m}  \E[L(k+l,\tau^\u_{m+l}],  &&
                                           -l < k < m
\end{align}
for $l$, $m \geq 1$
\qed
\end{lem}
We can use these occupation times to express the distribution of the levels visited upon exit from the upper boundary of an interval: by conditioning on the last level visited before exit, we obtain
\begin{equation}
   \label{e:exit}
\P[\tau^\u_m < \tau_{-l}, X_{\tau^\u_m} = m +u, J_{\tau^\u_m}]
= \sum_{1 \leq \nu \leq m+l-1} L(m-\nu, \tau_{-l} \wedge  \tau^\u_m)) A_{\nu+u}.
\end{equation}

If the transition matrix $A_{-1}$ is nonsingular, we may rewrite the
expressions in Lemma~\ref{t:nonnull}, for a null-recurrent as well as a
transient process $(X,J)$,
 in terms of the scale matrices $W(k) = \Xi(k)R^{-k}$ and prove the
 following, remembering that $W(k) = 0$ for $k \leq 0$.

\begin{lem}
   \label{t:Rnonsingular}
Assume that $A_{-1}$ is nonsingular. Irrespective of whether the process $(X,J)$ is transient or null-recurrent, for $l, m \geq 1$.
\begin{align}
   \label{e:un}
\E[L(k,\tau^\u_m)] & = W(m) R^{m-k} - W(k), &&
                                           k < m,
\\
   \label{e:deux}
\E[L(k, \tau_{-l}]  & = G^l W(k+l) - W(k,  &&
                                           k > -l,
\\
   \label{e:intervalB}
\E[L(k, \tau_{-l} \wedge  \tau^\u_m)]  & =
   W(m) W(m+l)^{-1} W(k+l) - W(k),  &&
                                           -l < k < m.
\end{align}
\end{lem}
\begin{proof}
For $k \geq 1$, 
\begin{align*}
\E[L(k,\tau^\u_m)] & = \sum_{0 \leq \nu \leq m-k-1} \Phi(k+\nu) R^\nu
                     \qquad \mbox{by (\ref{e:lim})}
  \\
 & = (\Xi(m)-\Xi(k)) R^{-k}
  \\
 & =  W(m) R^{m-k} - W(k)
\end{align*} 
by Lemma \ref{t:lattB}.  For $k \leq 0$,
\begin{align*}
\E[L(k,\tau^\u_m)] & = D_{-k,m} \Xi(m-k)  
  \\
& = \Xi(m) R^{-k}  \qquad \qquad \mbox{by Lemma \ref{t:lattA}}
  \\
 & = W(m) R^{m-k} - W(k)
\end{align*}
using Lemma \ref{t:lattB} again, and the fact that $W(k)=0$.  This completes the proof of (\ref{e:un}).

Next, for $k \geq 1$, conditioning on the first visit to level $-k$,
\[
\E[L(0, \tau_{-(k+l)}]  = \E[L(0, \tau_{-k})] + G^k \E[L(k, \tau_{-l})]
\]
or $\Theta(k+l) = \Theta(k) + G^k \E[L(k, \tau_{-l})]$ which, together with (\ref{e:thetaxi}) proves (\ref{e:deux}).   
For $-l \leq 0$, $\E[L(k, \tau_{-l})] = G^{-k} \E[L(0, \tau_{-l-k})]$ and (\ref{e:deux}) readily follows.

Finally, starting from
\[
\E[L(k, \tau^\u_m)]  = \E[L(k, \tau_{-l} \wedge  \tau^\u_m)] + D_{l,m} \E[L(k+l, \tau^\u_{m+l})],
\]
we obtain (\ref{e:intervalB}) after some simple manipulation.
\end{proof}

We observe that if we replace the scale function by its expression
$W(k)=\Xi(k) R^{-k}$ in the right-hand side of (\ref{e:deux}), we
find that  
\[
\E[L(-r, \tau_{-l}] R^{l-r}= G^l \E[L(0, \tau^\u_{l-r}], \qquad \mbox{ for $0 \leq r < l$.}
\]
This is a generalisation of (\ref{e:thetaxi}), both sides representing
the expected time in level $-l$, starting from level 0, after a first
passage to level $-l$ and before the next visit to level $-r$.

\subsection{Extrema --- Wiener-Hopf factorisation}

\paragraph{The lattice case.}

We assume in this section that the generator $Q$ is defective, and we
denote the killing time of the process by
\[
\zeta=\sup\{t\geq 0:J_t\neq \dagger\}.
\]
For simplicity of notation we
write $(X_\zeta,J_\zeta)$ for the value
$(X_{\zeta-},J_{\zeta-})$ of the process at the killing time  and we
define 
$\overline X = \sup_{t \leq \zeta} X_t$  and $\underline X=\inf_{t \leq \zeta} X_t$.

\begin{lem}
   \label{t:extrema}
For $m,l\geq 0$, 
\begin{align}
\p[\overline X=m,X_\zeta=m-l,J_\zeta]&= \E[L(m,\tau^\u_{m+1})]\, R^l \,\Delta_{\bs q}  
    \label{e:overlineX}\\ 
\p[\underline X=-m,X_\zeta=-m+l,J_\zeta]&=G^m \, \e[
       L(l,\tau_{-1})] \, \Delta_{\bs q}.
\label{e:underlineX}
\end{align}
\end{lem}
\begin{proof}
We have
\[
\p[\overline X=m,X_\zeta=m-l,J_\zeta] =\P[\tau_m = \tau^\u_m]\Phi(0)R^l\Delta_{\bs q}.
\]
Indeed, the first factor means that the process reaches level $m$
and does so before visiting any higher level, the product $\Phi(0) R^l
\Delta_{\bs q}$ is the probability that,  after getting to level $m$,  the process reaches the
cemetery state while being in level $m-l$ before having moved to any
level above $m$.  The product $\P[\tau_m = \tau^\u_m]\Phi(0)$ is equal
to $\Phi(m)$ by (\ref{e:phim})  and this proves the first statement.

The proof of (\ref{e:underlineX}) is similar: the process goes down
to level $-m$ with probability $G^m$, and $\e[L(l,\tau_{-1})] \,
\Delta_{\bs q} $ is the probability that it gets absorbed into the
cemetery state while visiting level $-m+l$ without having visited
level $-m-1$.
\end{proof}

The next statement is a direct consequence of Lemmas
\ref{t:Rnonsingular} and \ref{t:extrema}. 
\begin{cor}
   \label{t:extremabis}
If the matrix $A_{-1}$ is nonsingular, then 
\begin{align}
\p[\overline X=m,X_\zeta=m-l,J_\zeta]&=(W(m+1)R-W(m))R^l\Delta_{\bs q}, \label{e:overlineXbis}\\ 
\p[\underline X=-m,X_\zeta = -m+l ,J_\zeta] & =G^{m}(GW(l+1)-W(l))\Delta_{\bs q},
\label{e:underlineXbis}
\end{align}
for $m,l\geq 0$.
\qed
\end{cor}

\begin{rem} \em
Lemma \ref{t:extrema} is the generalisation to lattice MAPs of the Wiener-Hopf
factorisation of \levy processes (\cite[Page 165]{bertoin}).
\end{rem}

\section{Conclusion}
\label{s:conclusion}

We have presented the fundamental fluctuation theory concerning
one-sided MAPs in both lattice and non-lattice cases.  This work
demonstrates that it is very helpful to consider both theories in a
unified way, and to use ideas from one in thinking about the other.

We conclude with a final word about our assumptions.

\subsection{Jump-free and skip-free processes}\label{sec:MMBM}
It is well known that a MAP which is skip (jump) free in both
directions allows for a much simpler and more explicit analysis.  Such
a MAP is referred to as a \emph{quasi birth-death processes} (QBD) in
the lattice case, and as a \emph{Markov modulated Brownian motion}
(MMBM) or as a \emph{second order fluid model} in the non-lattice
case.  Various quantities of interest can be expressed in terms of the
fundamental matrices $G$, $R$ and their analogues $G_-$, $R_-$ for the
process $(-X,J)$ which is also skip (jump) free.
In particular, we have simple explicit formulas for~$H$.

\subsection{Negative jumps of phase type}
The framework of one-sided MAPs may be used to analyze more general
processes.  In this context \emph{phase type} (PH) distributions play
an important role.  A PH distribution is a
distribution of the life time of a 
transient Markov chain on finitely many states,
see~\cite[{\S}III.4]{APQ}, \cite[Chapter 2]{lara99}.  This property
allows us to represent MAPs 
with negative PH jumps as one-sided MAPs with an expanded set of
phases, the same idea is often used with respect to \levy processes
with PH jumps.  See for instance~\cite{asmussen_levy_PH,dbdlr05} 
and references therein.

\subsection{Subordinators in the non-lattice case}
   \label{s:subordinators}

   In the general case, where $X^{(i)}$ is allowed to be a
   subordinator for some phases $i$, one usually partitions the phases
   according to $E=E_\sim\cup E_\uparrow$, where the indices in
   $E_\uparrow$ correspond to subordinators, and we let $N_\sim$ and
   $N_\uparrow$ be the respective cardinalities. 

Considering the first passage Markov chain, we note that the level
$-x<0$ cannot be hit in any state $j\in E_\uparrow$  and so
\[
\P[J_{\tau_{-x}}] = \vligne{I \\ \Pi} \vligne{e^{\tilde G x} & 0},
\]
where $I$ is the $N_\sim \times N_\sim$ identity matrix, 
$\Pi$ is a transition probability matrix with dimensions $ N_\uparrow
\times N_\sim$ and $\tilde G$ is an irreducible intensity matrix  on
$E_\sim$. The pair $(\Pi, \tilde G)$ is identified by the matrix integral equation, 
\begin{align}
\nonumber
   \frac{1}{2}\Delta^2_{\bs\sigma}    
\vligne{I \\ \Pi} \tilde G^2 
& +\Delta_{\bs a}  \vligne{I \\ \Pi} \tilde G
+\int_0^\infty\Delta_{\bs \nu(\D  x)}  \vligne{I \\ \Pi}  \left(
  e^{\tilde Gx}-\matI- \tilde Gx\1{x<1}\right)
\\
   \label{eq:FG_gen}
& + \int_0^\infty (\tilde Q\circ U(\D   x)) \vligne{I \\ \Pi} e^{\tilde
  G x}= 0,
\end{align}
see~\cite[Thm. 2]{lambda}.
 
The rows of the matrix $R(x)$, $x>0$, corresponding to $E_\uparrow$ are 0,
because the process must hit level 0 in
some phase of $E_\sim$ on the way to the level $-x$ if it starts in $i\in
E_\uparrow$. 
Hence,
\[
e^{R x} = \vligne{e^{\tilde Rx} \\ 0}  \vligne{I & \Psi}
\]
where $\tilde R$ is an $N_\sim\times N_\sim$ matrix, $I$ is the
$N_\sim\times N_\sim$  identity matrix and $\Psi$ is a $N_\sim \times
N_\uparrow$ transition probability matrix.  The pair $(\tilde R,
\Psi)$  is a solution of
\begin{align}
\nonumber
      {\tilde R}^2 \vligne{I & \Psi}\frac{1}{2}\Delta^2_{\bs\sigma} +
 & {\tilde R}\vligne{I & \Psi}\Delta_{\bs
        a}+\int_0^\infty\left( e^{{\tilde R} x}-\matI-{\tilde R}x\1{x<1}\right)\vligne{I & \Psi}\Delta_{\bs \nu(\D
     x)}
\\
   \label{eq:FR_gen}
  & + \int_0^\infty e^{{\tilde R} x}\vligne{I & \Psi}(Q\circ U(\D
   x))=\matO.
\end{align}

\appendix
\section{The role of killing in Remark \ref{rem:4.2}}
\label{a:killing}
Consider the case where there is no value of $z$ (or $\a$) for which the convergence conditions of (\ref{eq:lattHseries}), respectively (\ref{eq:nonlattHint}) hold. In such a situation, it is tempting to introduce killing at a rate $q$ to transform the model into one in which we can apply the conclusions of Proposition \ref{prop:Hx}, and then let $q\to 0^+$.

For any killing parameter $q>0$, we can use the obvious physical justification to conclude that there is a family of matrices $H(k;q)$ such that, for any given $k$, $H(k;q)$ is a decreasing function of $q$ and 

\begin{equation}
   \label{e:hkq}
H(k;q) \rightarrow H(k) \quad \mbox{as $q \rightarrow 0^+$.} 
\end{equation}

It follows from convergence of eigenvalues and the
fact that $F_q(1)=Q-qI$ that there exist positive constants $\epsilon_q$ such
that the hypothesis of Proposition \ref{prop:Hx} on eigenvalues is satisfied for $z$ in some
interval $(1-\epsilon_q,1)$, and so
\begin{equation}
   \label{e:Hq}
H_q(z)  \equiv \sum_{k\in\mathbb Z} z^{k}H(k;q) =-z F_q(z)^{-1}
\end{equation}
for all $z \in (1-\epsilon_q,1)$.  

With $\DD$ the open unit disc in $\CC$, we can use the continuity with respect to $z$ of the eigenvalues of
$F_q(z)$ to extend the interval $(1-\epsilon_q,1)$ into an open
set $\sDD_q$ with
\[
(1-\epsilon_q, 1) \subset \sDD_q \subset \DD \setminus \{0\},
\]
such that for $z \in \sDD_q$, all eigenvalues of $F_q(z)/z$ have  negative real part.  We note that
$\sDD_q$ might converge to a non-empty set as $q \to 0^+$, but if so its limit does not contain any $z$ from $(0,1)$.

The proof of Proposition~4.2 does not depend on $z$ being real, and so
\begin{equation}
   \label{e:Hqbis}
H_q(z) = -z F_q(z)^{-1} \qquad \mbox{for $z \in \sDD_q$.}
\end{equation}

The functions
$F_q(z)$, $q \geq 0$ are analytic
in $\DD$
and their determinants are rational functions.  Therefore, for each $q$ there is
 a finite set $\mathcal Z_q$ of values of $z$ such that
$F_q(z)$ is singular, and the functions $z F_q(z)^{-1}$  are analytic
in $\DD \setminus \mathcal Z_q$.

Define $A_q(z)$ to be the analytic continuation of $H_q(z)$ in $\DD$.  We find  from \eqref{e:Hqbis} that
$
A_q(z) \equiv -z F_q(z)^{-1} 
$
for $z \in \DD \setminus \mathcal Z_q$. 
In other words, $A_q(z)$ is merely another name for the function
$-z F_q(z)^{-1}$ and \eqref{e:Hq} tells us that $-z F_q(z)^{-1}$ can be expressed
as the power series $H_q(z)$ over the part of its domain $\sDD_q$ where the series converges.

Using the physical interpretation (5) of $F(z)$, we know that
\begin{equation}
   \label{e:hqz}
\lim_{q\rightarrow 0}A_q(z) = -z F(z)^{-1} 
\end{equation}
for $z  \in \DD$ except where $F(z)$ is singular.  However if the hypothesis of Proposition \ref{prop:Hx} is not satisfied for any $z$, then $H(z)$ does not converge for any $z$ and so there is no basis to define an analytic continuation $A(z)$. It follows that we cannot think of \eqref{e:hqz} as saying that $\lim_{q\rightarrow 0}A_q(z) = A(z)$.  
In short, the introduction of
killing tells us something about the {\em individual}
$H(k)$s through~(\ref{e:hkq}), but does not extend to the generating function $H(z)$.

%

\paragraph{Acknowledgments}

The first author greatfully acknowledges the support of Sapere Aude Starting Grant 8049-00021B entitled “Distributional Robustness in Assessment of Extreme Risk. The third author would like to thank the Australian Research Council (ARC) for supporting this work through Laureate Fellowship FL130100039 and the ARC Centre of Excellence for Mathematical and Statistical Frontiers (ACEMS). These grants also supported visits to Melbourne by the other authors.

\bibliographystyle{abbrv}
\bibliography{ilt}

\end{document}